\documentclass[12pt,psamsfonts]{article}

\usepackage{theorem}
\usepackage{amssymb}
\usepackage{amsmath}
\usepackage{graphics}

\newtheorem{thm}{Theorem}[section]
\newtheorem{prop}[thm]{Proposition}
\newtheorem{cor}[thm]{Corollary}
\newtheorem{lm}[thm]{Lemma}
\newtheorem{hyp}[thm]{Hypotheses}
\newtheorem{conj}[thm]{Conjecture}

{\theorembodyfont{\rm}

\newtheorem{defn}[thm]{Definition}

 }

\def\e{\epsilon}

\def\l{\lambda}

\def\ni{\noindent}
\def\spec{\operatorname{Spec}}

\def\ps{\vspace{4pt}}

\def\sqr#1#2{{\vcenter{\vbox{\hrule height.#2pt \hbox{\vrule
width.#2pt height #1pt \kern #1pt \vrule width.#2pt}\hrule
height.#2pt}}}}

\def\qed{{\hfil \break  \rightline{$\sqr74$}}}

\def\C{{\mathbb C}}

\def\P{{\mathbb P}}

\def\cD{{\mathcal D}}

\def\cO{{\mathcal O}}

\def\cO{{\mathcal O}}

\def\e{\epsilon}

\def\mg{\overline M_{g,0}}

\begin{document}

\centerline{\bf  {\large Families of Rationally Connected Varieties}}

\

\

\centerline{\today}

\

\noindent {\bf Tom Graber}

\noindent Mathematics Department, Harvard University,

\noindent 1 Oxford st., Cambridge MA 02138, USA

\noindent graber{\char'100}math.harvard.edu

\

\noindent {\bf Joe Harris}

\noindent Mathematics Department, Harvard University,

\noindent 1 Oxford st., Cambridge MA 02138, USA

\noindent harris{\char'100}math.harvard.edu

\

\noindent {\bf Jason Starr}

\noindent Mathematics Department, M.I.T.,

\noindent Cambridge MA 02139, USA

\noindent jstarr{\char'100}math.mit.edu

\

\tableofcontents

\section{Introduction}

\subsection{Statement of results}

We will work throughout over the complex numbers, so that the results
here apply over any algebraically closed field of characteristic 0.

\ps

Recall that a smooth projective variety $X$ is said to be {\em
rationally connected} if two general points $p, q \in X$ can be joined
by a chain of rational curves. In case $\dim(X) \geq 3$ this turns out
to be equivalent to the a priori stronger condition that for any
finite subset $\Gamma \subset X$ there is a smooth rational curve $C
\subset X$ containing $\Gamma$ and having ample normal bundle.

\ps

Rationally connected varieties form an important class of
varieties. In dimensions 1 and 2 rational connectivity coincides with
rationality, but the two notions diverge in higher dimensions and in
virtually every respect the class of rationally connected varieties is
better behaved. For example, the condition of rational connectivity is
both open and closed in smooth proper families; there are geometric
criteria for rational connectivity (e.g. any smooth projective variety
with negative canonical bundle is rationally connected, so we know in
particular that a smooth hypersurface $X \subset \P^n$ of degree $d$
will be rationally connected if and only if $d \leq n$), and there
are, at least conjecturally, numerical criteria for rational
connectivity (see Conjecture~\ref{mumford} below).  In this paper we
will prove a conjecture of Koll\'ar, Miyaoka and Mori that represents
one more basic property of rational connectivity (also one not shared
by rationality): that if $X \to Y$ is a morphism with rationally
connected image and fibers, then the domain $X$ is rationally
connected as well. This will be a corollary of our main theorem:

\ps

\begin{thm}\label{mainth}
Let $f : X \to B$ be a morphism from a smooth projective variety to a
smooth projective curve over $\C$. If the general fiber of $f$ is
rationally connected, then $f$ has a section.
\end{thm}

Since this is really a statement about the birational equivalence
class of the morphism $f$, we can restate it in the equivalent form

\begin{thm}
If $K$ is the function field of a curve over $\C$, any rationally
connected variety $X$ defined over $K$ has a $K$-rational point.
\end{thm}

In this form, the theorem directly generalizes Tsen's theorem, which
is exactly this statement for $X$ a smooth hypersurface of degree $d
\leq n$ in projective space $\P^n$ (or more generally a smooth
complete intersection in projective space with negative canonical
bundle). It would be interesting to know if in fact rationally
connected varieties over other $C_1$ fields necessarily have rational
points.

\ps

As we indicated, one basic corollary of our main theorem is

\begin{cor}\label{totalspace}
Let $f : X \to Y$ be any dominant morphism of varieties. If $Y$ and
the general fiber of $f$ are rationally connected, then $X$ is
rationally connected.
\end{cor}

\ni {\em Proof}. We can assume (in characteristic 0, at least) that
$X$ and $Y$ are smooth projective varieties. Let $p$ and $q$ be
general points of $X$. We can find a smooth rational curve $C \subset
Y$ joining $f(p)$ and $f(q)$; let $X' = f^{-1}(C)$ be the inverse
image of $C$ in $X$. By Theorem~\ref{mainth}, there is a section $D$
of $X'$ over $C$. We can then connect $p$ to $q$ by a chain of
rational curves in $X'$ in three stages: connect $p$ to the point $D
\cap X_p$ of intersection of $D$ with the fiber $X_p$ of $f$ through
$p$ by a rational curve; connect $D \cap X_p$ to $D \cap X_q$ by $D$,
and connect $D \cap X_q$ to $q$ by a rational curve in $X_q$. \qed

There is a further corollary of Theorem~\ref{mainth} based on a
construction of Campana and Koll\'ar--Miyaoka--Mori: the {\it maximal
rationally connected fibration} associated to a variety $X$ (see [Ca],
[K] or [KMM]). Briefly, the maximal rationally connected fibration
associates to a variety $X$ a (birational isomorphism class of)
variety $Z$ and a rational map $\phi : X \to Z$ with the properties
that
\begin{itemize}
\item the fibers $X_z$ of $\phi$ are rationally connected; and conversely 
\item almost all the rational curves in $X$ lie in fibers of $\phi$:
for a very general point $z \in Z$ any rational curve in $X$ meeting
$X_z$ lies in $X_z$.
\end{itemize}
\ni The variety $Z$ and morphism $\phi$ are unique up to birational
isomorphism, and are called the \textit{mrc quotient} and \textit{mrc
fibration} of $X$, respectively.  They measure the failure of $X$ to
be rationally connected: if $X$ is rationally connected, $Z$ is a
point, while if $X$ is not uniruled we have $Z=X$. As observed in
Koll\'ar ([K], IV.5.6.3), we have the following Corollary:

\begin{cor}\label{quot}
Let $X$ be any variety and $\phi : X \to Z$ its maximal rationally
connected fibration. Then $Z$ is not uniruled.
\end{cor}

\ni {\em Proof}. Suppose that $Z$ were uniruled, so that through a
general point $z \in Z$ we could find a rational curve $C \subset Z$
through $z$. By Corollary~\ref{totalspace}, the inverse image
$\phi^{-1}(C)$ will be rationally connected, which means that every
point of the fiber $X_z$ will lie on a rational curve not contained in
$X_z$, contradicting the second defining property of mrc fibrations.
\qed

There are conjectured numerical criteria for a variety $X$ to be
either uniruled or rationally connected. They are

\begin{conj}\label{uniruled}
Let $X$ be a smooth projective variety. Then $X$ is uniruled if and
only if $H^0(X,K_X^m)=0$ for all $m > 0$.
\end{conj}

\ni and

\begin{conj}\label{mumford}
Let $X$ be a smooth projective variety. Then $X$ is rationally
connected if and only if $H^0(X,(\Omega^1_X)^{\otimes m})=0$ for all
$m > 0$.
\end{conj}

For each of these conjectures, the ``only if" part is known, and
straightforward to prove; the ``if" part represents a very difficult
open problem (see for example [K], IV.1.12 and IV.3.8.1). As another
consequence of our main theorem, we have an implication:

\begin{cor}
Conjecture~\ref{uniruled} implies Conjecture~\ref{mumford}
\end{cor}

\ni {\em Proof}. Let $X$ be any smooth projective variety that is not
rationally connected; assuming the statement of
Conjecture~\ref{uniruled}, we want to show that
$H^0(X,(\Omega^1_X)^{\otimes m})\neq 0$ for some $m > 0$. Let $\phi :
X \to Z$ be the mrc fibration of $X$. By hypothesis $Z$ has dimension
$n >0$, and by Corollary~\ref{quot} $Z$ is not uniruled. If we assume
Conjecture~\ref{uniruled}, then, we must have a non-zero section
$\sigma \in H^0(Z,K_Z^m)$ for some $m > 0$. But the line bundle
$K_Z^m$ is a summand of the tensor power $(\Omega^1_Z)^{\otimes nm}$,
so we can view $\sigma$ as a global section of that sheaf; pulling it
back via $\phi$, we get a nonzero global section of
$(\Omega^1_X)^{\otimes nm}$ \qed

\ps

\ni {\bf Acknowledgments}. We would like to thank Johan deJong,
 J\'anos Koll\'ar and Barry Mazur for many conversations, which were
 of tremendous help to us.

\ps  

\section{Preliminary definitions and constructions}

We will be dealing with morphisms $\pi : X \to B$ satisfying a number
of hypotheses, which we collect here for future reference. In
particular, for the bulk of this paper we will deal with the case $B
\cong \P^1$; we will show in section~\ref{barb} below both that the
statement for $B \cong \P^1$ implies the full Theorem~\ref{mainth}
and, as well, how to modify the argument that follows to apply to
general $B$.

\begin{hyp}\label{hyp}
$\pi : X \to B$ is a nonconstant morphism of smooth connected
projective varieties over $\C$, with $B \cong \P^1$. For general $b
\in B$, the fiber $X_b = \pi^{-1}(b)$ is rationally connected of
dimension at least 2.
\end{hyp}

Now suppose we have a class $\beta \in N_1(X)$ having intersection
number $d$ with a fiber of the map $\pi$.  We have then a natural
morphism
$$
\varphi : \mg(X,\beta) \to \mg(B,d)
$$
defined by composing a map $f : C \to X$ with $\pi$ and collapsing
components of $C$ as necessary to make the composition $\pi f$ stable.

\begin{defn}
Let $\pi : X \to B$ be a morphism satisfying \ref{hyp}, and let $f : C
\to X$ be a stable map from a nodal curve $C$ of genus $g$ to $X$ with
class $f_*[C] = \beta$. We say that $f$ is {\it flexible} relative to
$\pi$ if the map $\varphi : \mg(X,\beta) \to \mg(B,d)$ is dominant at
the point $[f] \in \mg(X,\beta)$; that is, if any neighborhood of
$[f]$ in $\mg(X,\beta)$ dominates a neighborhood of $[\pi f]$ in
$\mg(B,d)$.
\end{defn}

Now, it's a classical fact that the variety $\mg(B,d)$ has a unique
irreducible component whose general member corresponds to a flat map
$f : C \to B$ (see for example [C] and [H] for a proof). Since the map
$\varphi : \mg(X,\beta) \to \mg(B,d)$ is proper, it follows that if
$\pi : X \to B$ admits a flexible curve then $\varphi$ will be
surjective.  Moreover, $\mg(B,d)$ contains points $[f]$ corresponding
to maps $f : C \to B$ with the property that every irreducible
component of $C$ on which $f$ is nonconstant maps isomorphically via
$f$ to $B$. (For example, we could simply start with $d$ disjoint
copies $C_1,\dots,C_d$ of $B$ (with $f$ mapping each isomorphically to
$B$) and identify $d+g-1$ pairs of points on the $C_i$, each pair
lying over the same point of $B$.

\begin{prop}
If $\pi : X \to B$ is a morphism satisfying \ref{hyp} and $f : C \to
X$ a flexible stable map, then $\pi$ has a section.
\end{prop}

Our goal in what follows, accordingly, will be to construct a flexible
curve $f : C \to X$ for an arbitrary $\pi : X \to B$ satisfying
\ref{hyp}.

\subsection{The first construction}

To manufacture our flexible curve, we apply two basic constructions,
which we describe here. (These constructions, especially the first,
are pretty standard: see for example section II.7 of [K].)  We start
with a basic lemma:

\begin{lm}\label{bundle}
Let $C$ be a smooth curve and $E$ any vector bundle on $C$; let $n$ be
any positive integer. Let $p_1,\dots,p_N \in C$ be general points and
$\xi_i \subset E_{p_i}$ a general one-dimensional subspace of the
fiber of $E$ at $p_i$; let $E'$ be the sheaf of rational sections of
$E$ having at most a simple pole at $p_i$ in the direction $\xi_i$ and
regular elsewhere. For $N$ sufficiently large we will have
$$
H^1(C, E'(-q_1-\dots-q_n)) = 0
$$ 
for any $n$ points $q_1,\dots,q_n \in C$.
\end{lm}

\ni {\it Proof}. To start with, we will prove simply that
$H^1(C,E')=0$. Since this is an open condition, it will suffice to
exhibit a particular choice of points $p_i$ and subspaces $\xi_i$ that
works.  Denoting the rank of $E$ by $r$, we take $N = mr$ divisible by
$r$ and choose $m$ points $t_1,\dots,t_m \in C$. We then specialize to
the case
\begin{align*}
p_1 = \dots = p_r &= t_1; & \xi_1,\dots,\xi_r &\mbox{ spanning }
E_{t_1} \\ 
p_{r+1} = \dots = p_{2r} &= t_2; & \xi_{r+1},\dots,\xi_{2r}
&\mbox{ spanning } E_{t_2}
\end{align*}
and so on. In this case we have $E' = E(t_1 + \dots + t_m)$, which we
know has vanishing higher cohomology for sufficiently large $m$.

Given this, the statement of the lemma follows: to begin with, choose
any $g+n$ points $r_1,\dots,r_{g+n} \in C$. Applying the argument thus
far to the bundle $E(-r_1-\dots -r_{g+n})$, we find that for $N$
sufficiently large we will have $H^1(C, E'(-r_1-\dots-r_{g+n})) =
0$. But now for any points $q_1,\dots,q_n \in C$ we have
$$
q_1 + \dots + q_n = r_1+\dots + r_{g+n} - D
$$
for some effective divisor $D$ on $C$. It follows then that
\begin{align*}
h^1(C, E'(-q_1-\dots-q_n)) &= h^1(C, E'(-r_1-\dots -r_{g+n})(D)) \\
&\leq h^1(C, E'(-r_1-\dots -r_{g+n})) \\
&=0
\end{align*}
\qed

The relevance of this to our present circumstances will perhaps be
made clear by the following:

\begin{lm}\label{normal}
Let $X$ be a smooth projective variety, $C$ and $C'\subset X$ two
nodal curves meeting at points $p_1,\dots,p_\delta$; suppose $C$ and
$C'$ are smooth with distinct tangent lines at each point $p_i$. Let
$D = C \cup C'$ be the union of $C$ and $C'$; and let $N_{C/X}$ and
$N_{D/X}$ be the normal sheaves of $C$ and $D$ in $X$. We have then an
inclusion of sheaves
$$
0 \to N_{C/X} \to N_{D/X}|_C
$$
identifying the sheaf of sections of $N_{D/X}|_C$ with the sheaf of
rational sections of $N_{C/X}$ having at most a simple pole at $p_i$
in the normal direction determined by $T_{p_i}C'$. Moreover, if
$\tilde D \subset \spec \C[\e]/(\e^2) \times X$ is a first-order
deformation of $D$ in $X$ corresponding to a global section $\sigma
\in H^0(N_{D/X})$, then $\tilde D$ smooths the node of $D$ at $p_i$ if
and only if the restriction $\sigma|_U$ of $\sigma$ to a neighborhood
$U$ of $p$ in $C$ is not in the image of $N_{C/X}$.
\end{lm}

\ps

Now suppose $\pi : X \to B$ is a morphism satisfying our basic
hypotheses~\ref{hyp}, and $C \subset X$ a smooth, irreducible curve of
genus $g$. For a general point $p \in C$, let $X_p = \pi^{-1}(\pi(p))$
be the fiber of $\pi$ through $p$. By hypothesis, $X_p$ is a smooth,
rationally connected variety, so that we can find a smooth rational
curve $C' \subset X_p$ meeting $C$ at $p$ (and nowhere else) with
arbitrarily specified tangent line at $p$, and having ample normal
bundle $N_{C'/X}$.

Choose a large number of general points $p_1,\dots,p_\delta \in C$,
and for each $i$ let $C_i \subset X_{p_i}$ be such a smooth rational
curve, with $T_{p_i}C_i$ a general tangent line to $X_{p_i}$ at
$p_i$. Combining the preceding two lemmas, we see that for $\delta$
sufficiently large, the normal bundle $N_{C'/X}$ of the union $C' = C
\cup (\cup C_i)$ will be generated by its global sections; in
particular, by Lemma~\ref{normal} there will be a smooth deformation
$\tilde C$ of $C'$. Moreover, for any given $n$ we can choose the
number $\delta$ large enough to ensure that $H^1(C,
N_{C'/X}|_C(-r_1-\dots-r_{g+n})) = 0$ for some $g + n$ points
$r_1,\dots,r_{g+n} \in C$; it follows that $H^1(\tilde C, N_{\tilde
C/X}(-r_1-\dots-r_{g+n})) = 0$ for some $r_1,\dots,r_{g+n} \in \tilde
C$ and hence that
$$
H^1(\tilde C, N_{\tilde C/X}(-q_1-\dots-q_n)) = 0
$$
for any $n$ points on $\tilde C$.

The process of taking a curve $C \subset X$, attaching rational curves
in fibers and smoothing to get a new curve $\tilde C$, is our first
construction. It has the properties that

\begin{enumerate}
\item the genus $g$ of the new curve $\tilde C$ is the same as the
genus of the curve $C$ we started with;
\item the degree $d$ of $\tilde C$ over $B$ is the same as the degree
of $C$ over $B$;
\item the branch divisor of the composite map $\tilde C
\hookrightarrow X \to B$ is a small deformation of the branch divisor
of $C \hookrightarrow X \to B$; and again,
\item for any $n$ points $q_1,\dots,q_n \in \tilde C$ we have
$H^1(\tilde C, N_{\tilde C/X}(-q_1-\dots-q_n)) = 0$
\end{enumerate}

Here is one application of this construction. Suppose we have a smooth
curve $C \subset X$ such that the projection $\mu : \pi|_C : C \to B$
is simply branched---that is, the branch divisor of $\mu$ consists of
$2d+2g-2$ distinct points in $B$---and such that each ramification
point $p \in C$ of $\mu$ is a smooth point of the fiber
$X_p$. Applying our first construction with $n = 2d+2g-2$, we arrive
at another smooth curve $\tilde C$ that is again simply branched over
$B$, with all ramification occurring at smooth points of fibers of
$\pi$. But now the condition that $H^1(\tilde C, N_{\tilde
C/X}(-q_1-\dots-q_n)) = 0$ applied to the $n = 2d+2g-2$ ramification
points of the map $\tilde \mu : \tilde C \to B$ says that if we pick a
normal vector $v_i$ to $\tilde C$ at each ramification point $p_i$ of
$\tilde \mu$ we can find a global section of the normal bundle
$N_{\tilde C/X}$ with value $v_i$ at $p_i$. Moreover, since
ramification occurs at smooth points of fibers of $\pi$, for any
tangent vectors $w_i$ to $B$ at the image points $\pi(p_i)$ we can
find tangent vectors $v_i \in T_{p_i}X$ with $d\pi(v_i) = w_i$. It
follows that {\it as we deform the curve $\tilde C$ in $X$, the branch
points of $\tilde \mu$ move independently}. A general deformation of
$\tilde C \subset X$ thus yields a general deformation of $\tilde
\mu$---in other words, the curve $\tilde C$ is flexible. We thus make
the

\begin{defn}
Let $\pi : X \to B$ be as in~\ref{hyp}, and let $C \subset X$ be a
smooth curve such that the projection $\mu : \pi|_C : C \to B$ is
simply branched. If each ramification point $p \in C$ of $\mu$ is a
smooth point of the fiber $X_p$ containing it, we will say the curve
$C$ is {\it pre-flexible}.
\end{defn}

\ni In these terms, we have established the

\begin{lm}\label{preflex}
Let $\pi : X \to B$ be as in~\ref{hyp}. If $X$ admits a pre-flexible
curve, the map $\pi$ has a section.
\end{lm}

\ni \underbar{Remark}. Note that we can extend the notion of
pre-flexible and the statement of Lemma~\ref{preflex} to stable maps
$f : C \to X$: we say that such a map is preflexible is the
composition $\pi f$ is simply branched, and for each ramification
point $p$ of $\pi f$ the image $f(p)$ is a smooth point of the map
$\pi$, the statement of Lemma~\ref{preflex} holds.

\subsection{The second construction}

Our second construction is a very minor modification of the
first. Given a family $\pi : X \to B$ as in~\ref{hyp} and a smooth
curve $C \subset X$, we pick a general fiber $X_b$ of $\pi$ and two
points $p, q \in C \cap X_b$. We then pick a rational curve $C_0
\subset X_b$ with ample normal bundle in $X_b$, passing through $p$
and $q$ and not meeting $C$ elsewhere. We also pick a large number $N$
of other general points $p_i \in C$ and rational curves $C_i \subset
X_{p_i}$ in the corresponding fibers, meeting $C$ just at $p_i$ and
having general tangent line at $p_i$. Finally, we let $C' = C \cup C_0
\cup (\cup C_i)$ be the union, and $\tilde C$ a smooth deformation of
$C'$ (as before, if we pick $N$ large enough, the normal bundle
$N_{C'/X}$ will be generated by global sections, so smoothings will
exist). This process, starting with the curve $C \subset X$ and
arriving at the new curve $\tilde C$, is our second construction. It
has the properties that

\begin{enumerate}
\item the degree $d$ of $\tilde C$ over $B$ is the same as the degree
of $C$ over $B$;
\item the genus of the new curve $\tilde C$ is one greater than the
genus of the curve $C$ we started with;
\item for any $n$ points $q_1,\dots,q_n \in \tilde C$ we have
$H^1(\tilde C, N_{\tilde C/X}(-q_1-\dots-q_n)) = 0$; and
\item the branch divisor of the composite map $\tilde C
\hookrightarrow X \to B$ has two new points: it consists of a small
deformation of the branch divisor of $C \hookrightarrow X \to B$,
together with a pair of simple branch points $b', b'' \in B$ near $b$,
each having as monodromy the transposition exchanging the sheets of
$\tilde C$ near $p$ and $q$.
\end{enumerate}

In effect, we have simply introduced two new simple branch points to
the cover $C \to B$, with assigned (though necessarily equal)
monodromy. Note that we can apply this construction repeatedly, to
introduce any number of (pairs of) additional branch points with
assigned (simple) monodromy; or we could carry out a more general
construction with a number of curves $C_0$.

\section{Proof of the main theorem}\label{mainproof}

\subsection{The proof in case $B = \P^1$}\label{mainarg}

We are now more than amply equipped to prove the theorem. We start
with a morphism $\pi : X \to B$ as in~\ref{hyp}. To begin with, by
hypothesis $X$ is projective; embed in a projective space and take the
intersection with $\dim(X) - 1$ general hyperplanes to arrive at a
smooth curve $C \subset X$. This is the curve we will start with.

What do $C$ and the associated map $\mu : C \hookrightarrow X \to B$
look like? To answer this, start with the simplest case: suppose that
the fibers $X_b$ of $\pi$ do not have multiple components, or in other
words that the singular locus $\pi_{\rm sing}$ of the map $\pi$ has
codimension 2 in $X$. In this case we are done: $C$ will miss
$\pi_{\rm sing}$ altogether, so that all ramification of $\mu : C \to
B$ will occur at smooth points of fibers; and simple dimension counts
show that the branching will be simple. In other words, $C$ will be
pre-flexible already.

The problems start if $\pi$ has multiple components of fibers. If $Z
\subset X_b$ is such a component, then each point $p \in C \cap Z$
will be a ramification point of $\mu$, and no deformation of $C$ will
move the corresponding branch point $\pi(p) \in B$. The curve $C$ can
not be flexible. And of course it's worse if $\pi$ has a multiple
(that is, everywhere-nonreduced) fiber: in that case $\pi$ cannot
possibly have a section.

\ps

To keep track of such points, let $M \subset B$ be the locus of points
such that the fiber $X_b$ has a multiple component. Outside of $M$,
the map $\mu : C \to B$ is simply branched, and all ramification
occurs at smooth points of fibers of $\pi$.

\ps

Now here is what we're going to do.  First, pick a base point $p_0 \in
B$, and draw a cut system: that is, a collection of real arcs joining
$p_0$ to the branch points $M \cup N$ of $\mu$, disjoint except at
$p_0$. The inverse image in $C$ of the complement $U$ of these arcs is
simply $d$ disjoint copies of $U$; call the set of sheets $\Gamma$
(or, if you prefer, label them with the integers 1 through $d$). Now,
for each point $b \in M$, denote the monodromy around the point $b$ by
$\sigma_b$, and express this permutation of $\Gamma$ as a product of
transpositions:
$$
\sigma_b = \tau_{b,1}\tau_{b,2}\dots\tau_{b,{k_b}} 
$$
so that in other words
$$
\tau_{b,{k_b}}\dots\tau_{b,2}\tau_{b,1}\sigma_b =  I
$$
is the identity. For future reference, let $k = \sum k_b$. We will
proceed in three stages.

\ps

\ni \underbar{Stage 1}: We use our second construction to produce a
new curve $\tilde C$ with, in a neighborhood of each point $b \in M$,
$k_b$ new pairs of simple branch points $s_{b,i}, t_{b,i} \in B$, with
the monodromy around $s_{b,i}$ and $t_{b,i}$ equal to
$\tau_{b,i}$. Note that $\tilde C$ will have genus $g(C) + k$, and
that the branch divisor of the projection $\tilde \mu : \tilde C \to
B$ will be the union of a deformation $\tilde N$ of $N$, the points
$s_{b,i}$ and $t_{b,i}$, and $M$. In particular we can find disjoint
discs $\Delta_b \subset B$, with $\Delta_b$ containing the points $b$
and $t_{b,1}, t_{b,2},\dots,t_{b,{k_b}}$, so that the monodromy around
the boundary $\partial \Delta_b$ of $\Delta_b$ is trivial.

\ps

Now, for any fixed integer $n$ this construction can be carried out so
that the curve $\tilde C$ has the property that $H^1(\tilde C,
N_{\tilde C/X}(-q_1-\dots-q_n)) = 0$ for any $n$ points $q_i \in
\tilde C$. Here we want to choose
$$
n = \#N +  2k
$$
so that there are global sections of the normal bundle $N_{\tilde
C/X}$ with arbitrarily assigned values on the ramification points of
$\tilde C$ over $N$ and the points $s_{b,i}$ and $t_{b,i}$.  This
means in particular that we can deform the curve $\tilde C$ so as to
deform the branch points of $\tilde \mu$ outside of $M$
independently. What we will do, then, is

\newpage

\ni \underbar{Stage 2}: We will vary $\tilde C$ so as to {\em keep all
the branch points $b \in N$ and all the points $s_{b,i}$ fixed; and
for each $b \in M$ specialize all the branch points $t_{b,i}$ to $b$
within the disc $\Delta_b$}.

%\vspace{3.5in}
\label{splat}
\begin{picture}(330,270)
\put(10,10){\makebox(330,240){\includegraphics{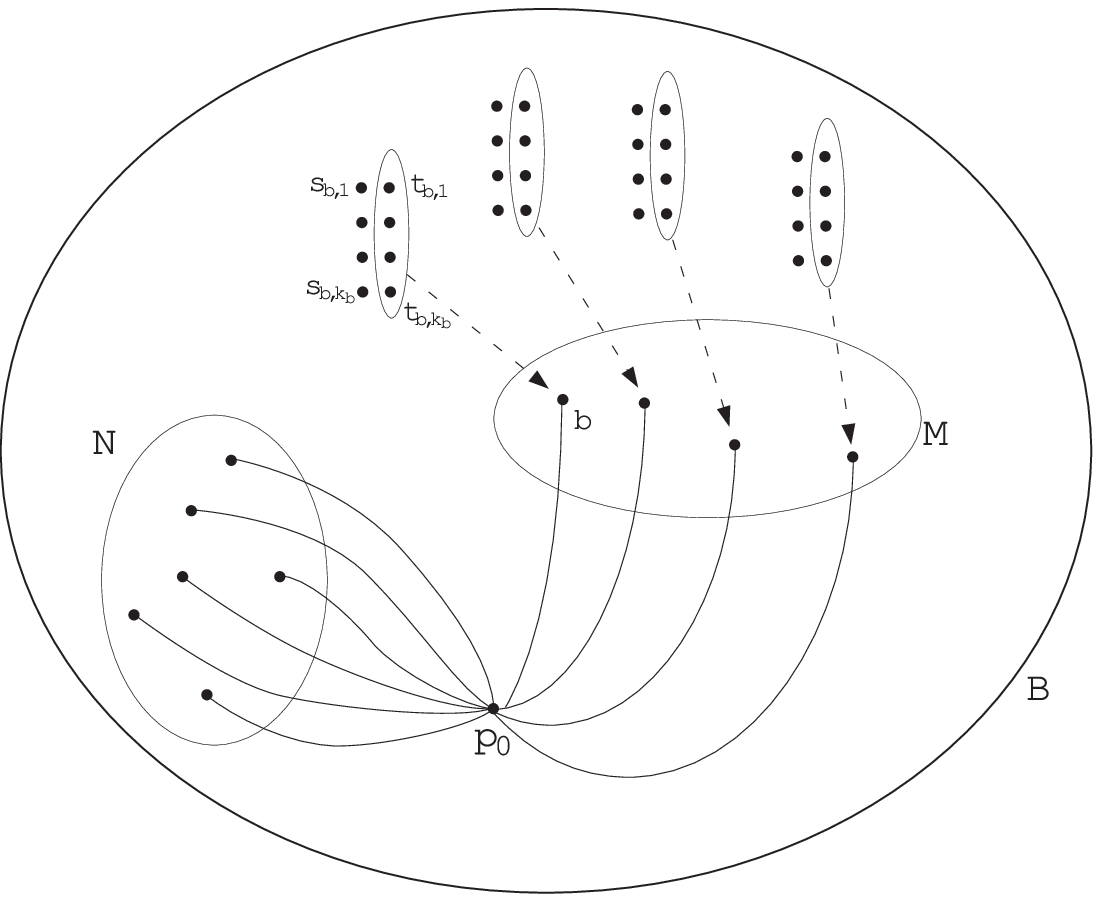}}}
\end{picture}
%\hspace{.2in}
%\special{picture splat}

\

\ni To say this more precisely, let $\beta \in N_1(X)$ be the class of
the curve $\tilde C$, and consider the maps
$$
\overline M_{g',0}(X, \beta) \longrightarrow \overline M_{g',0}(B, d)  \longrightarrow  B_{2d+2g'-2}
$$
with the second map assigning to a stable map $C \to B$ its branch
divisor.  What we are saying is, starting at the branch divisor
$$
D_1 =  \tilde N + \sum s_{b,i} +  \sum t_{b,i} + \sum_{b \in M} k_b\cdot b
$$
of the map $\tilde \mu$, draw an analytic arc $\gamma = \{D_\l\}$ in
the subvariety
$$
\Phi \, = \, \tilde N + \sum s_{b,i} + \sum_{b \in M} k_b\cdot b +
\sum(\Delta_b)_{k_b} \, \subset \, B_{2d+2g'-2}
$$
tending to the point
$$
D_0 = \tilde N + \sum s_{b,i} +  2\sum_{b \in M} k_b\cdot b .
$$ 
Since the image of the composition
$$
\overline M_{g',0}(X, \beta)  \longrightarrow  B_{2d+2g'-2}
$$
contains $\Phi$, we can find an arc $\delta = \{f_{\nu}\}$ in
$\overline M_{g',0}(X, \beta)$ that maps onto $\gamma$, with $f_1$ the
inclusion $\tilde C \hookrightarrow X$.

\ps

\ni \underbar{Stage 3}: Let $f_0 : C_0 \to X$ be the limit, in
$\overline M_{g',0}(X, \beta)$, of the family of curves constructed in
Stage 2: that is, the point of the arc $\delta$ over $D_0 \in \Phi
\subset B_{2d+2g'-2}$. Let $A \subset C_0$ be the normalization of any
irreducible component of $C_0$ on which the composition $\pi f_0$ is
nonconstant (that is, whose image is not contained in a fiber), and
let $f : A \to X$ be the restriction of $f_0$ to $A$.

\ps

By construction, the composition $\pi f$ is unramified over a
neighborhood of $M$: the monodromy around the boundary $\partial
\Delta_b$ of each disc $\Delta_b$ is trivial, and it can be branched
over at most one point $b$ inside $\Delta_b$, so it can't be branched
at all over $\Delta_b$. Indeed, it is (at most) simply branched over
each point of $N$ and each point $s_{b,i}$, and unramified
elsewhere. Moreover, since we can carry out the specialization of
$\tilde C$ above with the entire fiber of $\tilde C$ over the points
of $N$ and the $s_{b,i}$ fixed, the ramification of $\pi f$ on $A$
over these points will occur at smooth points of the corresponding
fibers of $\pi$. In other words, {\em the map $f : A \to X$ is
preflexible}, and we are done.  \qed

\subsection{The proof for arbitrary curves $B$}\label{barb}

As we indicated at the outset, there are two straightforward ways of
extending this result to the case of arbitrary curves $B$.

\ps

For one thing, virtually all of the argument we have made goes over
without change to the case of base curves $B$ of any genus $h$. The
one exception to this is the statement that the space $\overline
M_{g,0}(B,d)$ of stable maps $f : C \to B$ of degree $d$ from curves
$C$ of genus $g$ to $B$ has a unique irreducible component whose
general member corresponds to a flat map $f : C \to B$. This is false
in general---consider for example the case $g = d(h-1) + 1$ of
unramified covers. It is true, however, if we restrict ourselves to
the case $g \gg h,d$ (that is, we have a large number of branch
points) and look only at covers whose monodromy is the full symmetric
group $S_d$. Given this fact, and observing that our second
construction allows us to increase the number of branch points of our
covers $C \to B$ arbitrarily, the theorem can be proved for general
$B$ just as it is proved above for $B \cong \P^1$.

\ps

Alternatively, Johan deJong showed us a simple way to deduce the
theorem for general $B$ from the case $B \cong \P^1$ alone. We argue
as follows: given a map $\pi : X \to B$ with rationally connected
general fiber, we choose any map $g : B \to \P^1$ expressing $B$ as a
branched cover of $\P^1$. We can then form the ``norm" of $X$: this is
the (birational isomorphism class of) variety $Y \to \P^1$ whose fiber
over a general point $p \in \P^1$ is the product
$$
Y_p = \prod_{q \in g^{-1}(p)} X_q .
$$
Since the product of rationally connected varieties is again
rationally connected, it follows from the $\P^1$ case of the theorem
that $Y \to \P^1$ has a rational section, and hence so does $\pi$.

\section{An example}

There are a number of disquieting aspects of the argument in
Section~\ref{mainarg}, and in particular about the specialization in
Stage 2 of that argument. Clearly the curve $f : A \to X$ constructed
there cannot meet any multiple component of a fiber of $\pi : X \to
B$; that is, for each $b \in M$ it must meet the fiber $X_b$ only in
reduced components of $X_b$. This raises a number of questions: what
if the fiber $X_b$ is multiple?  How can the curve $\tilde C$, which
meets all the multiple components of $X_b$, specialize to one that
misses them all? And can we say which reduced components of $X_b$ the
curve $A$ will meet?

\ps

The answers to the first two questions are straightforward: in fact,
the argument given here proved that {\em the map $\pi : X \to B$
cannot have multiple fibers}, that is, every fiber $X_b$ must have a
reduced component.\footnote{In fact, this assertion is nearly
tantamount to Theorem~\ref{mainth} itself, but we were unable to prove
it directly except under additional and restrictive hypotheses} As for
the second, what must happen is that as our parameter $\delta \to 0$,
the points of intersection of $C_\delta$ with the multiple components
of $X_b$ slide toward the reduced components of $X_b$; the curve $C_0$
produced in the limit will have components contained in the fiber
$X_b$ and joining the points of intersection of $A$ with $X_b$ to each
of the multiple components. Finally, the answer to the third
question---and indeed the whole process---may be illuminated by
looking at a simple example; we will do this now.

\ps

To start, we have to find an example of a map $\pi : X \to B$ with
rationally connected general fiber and a special fiber having a
multiple component (and smooth total space $X$). Without question, the
simplest example will have general fiber $X_b \cong \P^1$, and special
fiber a chain of three smooth rational curves:

%\vspace{.3in}
\label{fiber}
\begin{picture}(330,30)
\put(10,0){\makebox(330,25){\includegraphics{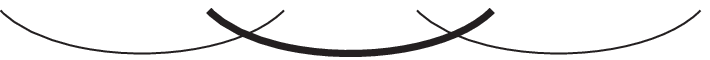}}}
\end{picture}
%\hspace{.8 in} \special{picture fiber}

\

The middle component will have multiplicity 2 in the fiber, and
self-intersection $-1$; the outer two components will each appear with
multiplicity 1 in the fiber, and will have self-intersection $-2$.
The simplest way to construct a family with such a fiber is to start
with a trivial family $X_0 = \P^1 \times \P^1 \to \P^1$, blow up any
point $p$, and then blow up the point $q$ of intersection of the
exceptional divisor with the proper transform of the fiber through $p$
to obtain $X$.

%\vspace{4.3in}
\begin{picture}(330,320)
\put(10,0){\makebox(330,310){\includegraphics{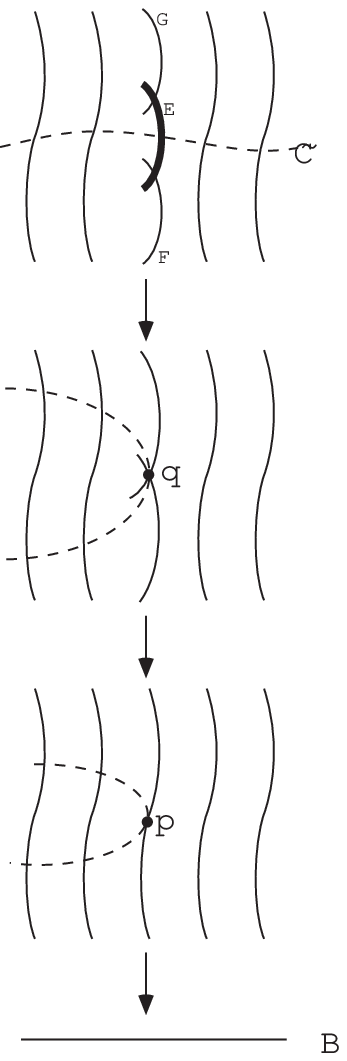}}}
\end{picture}
%\hspace{1.5in}\special{picture family}

\

We will denote by $F$ the proper transform in $X$ of the fiber through
$p$ in $\P^1 \times \P^1$, and by $G$ the proper transform of the
first exceptional divisor; the second exceptional divisor---the
multiple component of the special fiber---we will call $E$. To arrive
at the simplest possible curve $\tilde C \subset X$ meeting the
multiple component $E$ of the special fiber of this family, we start
with a curve $C \subset \P^1 \times \P^1$ of degree $2$ over $B$ that
is simply tangent to the special fiber at the point $p$; the proper
transform $\tilde C$ of $C$ in $X$ will then meet $E$ once
transversely and $F$ and $G$ not at all. (We're not trying to make
excuses here, but note that it's virtually impossible to draw a decent
picture of the configuration $\tilde C \subset X \to B$: the curve
$\tilde C$ is supposed to meet $E$ once transversely, but still have
degree 2 over $B$ and be ramified over $B$ at its point of
intersection with $E$.)

\ps

Now that we've got this set up, what happens when we push another
branch point of $\tilde C \to B$ in to the special fiber of $\pi$? The
answer is that one of three things can happen, two generically. We
will describe these first geometrically in terms of the original curve
$C \subset \P^1 \times \P^1$ and its proper transforms, and then write
down typical equations.

\ps

One possibility is that the ramification point $p$ of $C$ over $B$
becomes a node. In this case the limit $C_0$ of the proper transforms
$C_\nu$ of the curves will actually contain the component $E$ of the
special fiber (the limit of the proper transforms is not the proper
transform of the limiting curve, but rather its total transform minus
the divisor $E + G$). The remaining component---the actual proper
transform of the limiting curve---will have two distinct sheets in a
neighborhood of the special fiber, meeting $G$ transversely in
distinct points, and each of course unramified over $B$:

\newpage

\

%\vspace{4.3in}
\label{firstfam}
\begin{picture}(330,280)
\put(10,30){\makebox(330,250){\includegraphics{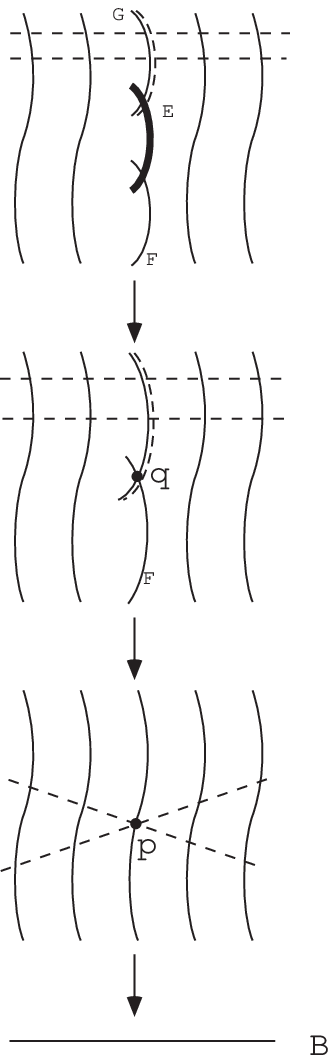}}}
\end{picture}
%\hspace{1.5in}\special{picture firstfam}

This specialization is easy to see in terms of equations: if we choose
affine coordinates $x$ on our base $\P^1$ and $y$ on the fiber, we can
write the equation of our family $\{C_\nu\}$ of curves as
$$
C_\nu \; : \; y^2 = x^2 - \nu x
$$ 
and specialize the branch point over $x = \nu$ simply by letting $\nu
\to 0$. We can see either from this family of equations, or
geometrically, that as $\nu$ tends to 0 the point of intersection of
the proper transform $\tilde C_\nu$ of $C_\nu$ slides along $E$ toward
the point of intersection $E \cap G$; when it reaches $E \cap G$ the
limiting curve becomes reducible, splitting off a copy of $G$.

\ps

Now, by the symmetry of $X \to B$---we could also blow down the curves
$E$ and $F$ in $X$ to obtain $\P^1 \times \P^1$---we would expect that
there would be a similar specialization with the roles of $F$ and $G$
reversed, and there is: if the curve $C \subset \P^1 \times \P^1$
specializes to one containing the fiber $\{0\} \times \P^1$, the limit
$\tilde C_0$ of the proper transforms will (generically) consist of
the union of $F$ with a curve $A$, with $A$ unramified of degree 2
over $B$ in a neighborhood of the special fiber and meeting the
special fiber in two distinct points of $F$.

\

%\vspace{4.3in}
\label{secondfam}
\begin{picture}(330,280)
\put(10,20){\makebox(330,260){\includegraphics{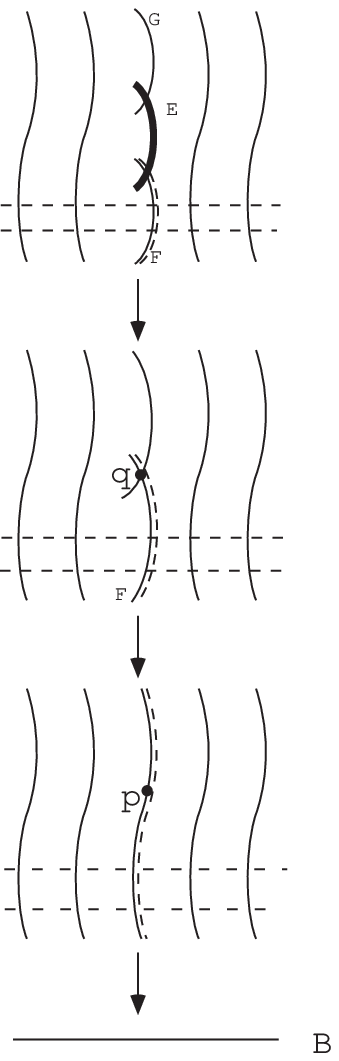}}}
\end{picture}
%\hspace{1.5in}\special{picture secondfam}

Finally, there is a common specialization of these two families: if
the curve $C \subset \P^1 \times \P^1$ specializes to one that both
contains the fiber $\{0\} \times \P^1$ and is singular at the point
$q$---that is, consists in a neighborhood of the special fiber of the
fiber and two sections, one passing through $p$---then the limit
$\tilde C_0$ of the proper transforms will consist of the union of all
three components $E$, $F$ and $G$ of the special fiber with a curve
$A$ consisting of a section meeting the special fiber in a point of
$F$ and a section meeting the special fiber in a point of $E$:

\newpage

\

%\vspace{4.3in}
\label{thirdfam}
\begin{picture}(330,280)
\put(10,20){\makebox(330,260){\includegraphics{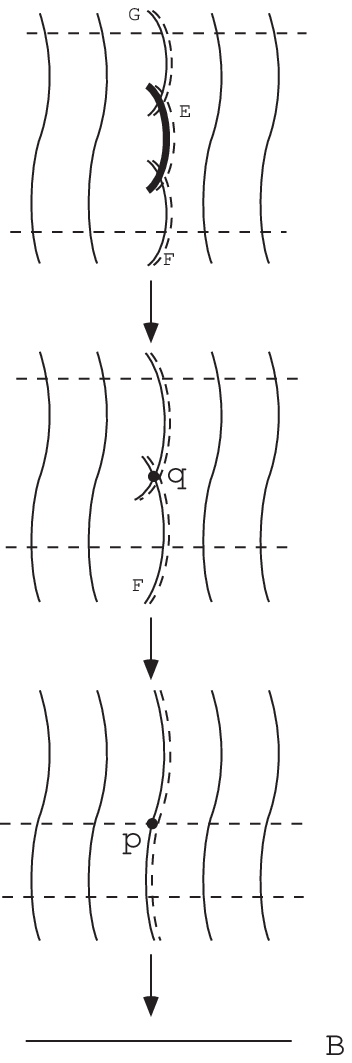}}}
\end{picture}
%\hspace{1.5in}\special{picture thirdfam}

It's also very instructive to look at this example from the point of
view of the equations of the curves. To begin with, denote by
$|\cO_X(d,e)|$ the total transform of the linear system of curves of
bidegree $(d,e)$ on $\P^1 \times \P^1$. We are looking here at the
linear system
$$
\tilde \cD = |\cO_X(1,2)(-G-2E)|,
$$
 that is, the proper transform of the linear series $\cD$ of curves $C
\subset \P^1 \times \P^1$ of bidegree $(1,2)$ that pass through $p$
with vertical tangent. Explicitly, these curves form a $3$-dimensional
linear series, which we may write in affine coordinates $(x,y)$ on
$\P^1 \times \P^1$ as
$$
\cD \, = \, \{a xy^2 + b y^2 + c xy + d x\}_{[a,b,c,d] \in \P^3}
$$
Writing the equation of a typical member of $\cD$ as a polynomial in
$y$:
$$
(a x + b) \cdot y^2 + (c x)\cdot y + (d x) = 0
$$
we see that its branch divisor is the zero locus of the quadratic polynomial
$$
(d x)^2 - 4(a x + b)(d x) \, = \, (c^2 - 4ad)\cdot x^2 - 4db\cdot x
$$
whose roots are at $x = 0$ and $x = 4db/(c^2 - 4ad)$. It's probably
best to express this in terms of the maps
$$
\overline M_{0,0}(X, [\tilde\cD]) \longrightarrow \overline M_{0,0}(B,
2) \longrightarrow B_2
$$
introduced in section~\ref{mainproof} above. Here, the variety
$\overline M_{0,0}(X, [\tilde\cD])$ has a component $M$ which is a
blow-up of the $\P^3_{[a,b,c,d]}$ parametrizing the linear series
$\tilde \cD$ (it also has a second, extraneous component whose general
point corresponds to a map $f : C \to X$ with reducible domain and
image containing the line $y=0$ doubly; this component is not involved
here). The image of the composite map $M \to B_2$ is simply the locus
$B_0 \cong B \subset B_2$ of divisors of degree 2 in $B \cong \P^1$
containing the point $x=0$, with the map
$$
\eta : M \to \P^3_{[a,b,c,d]} \to B_0 \cong \P^1
$$
given by
$$
[a,b,c,d] \mapsto [4db, c^2 - 4ad].
$$
What we see in particular from this is that {\em the fiber of $\eta$
over the point $x=0$ is reducible}, with components given by $d=0$ and
$b=0$. Now, in Stage 2 of our argument, as applied here, we start with
an arc $\gamma \subset B_0 \subset B_2$ in which the second branch
point approaches $x=0$, and lift that to an arc $\delta \subset M$.
If our arc $\delta \subset M$ lifting the arc $\gamma \subset B_0
\subset B_2$ approaches the component $d=0$---whose general member
corresponds to a curve $C \subset \P^1 \times \P^1$ singular at
$p$---we get a family of stable maps whose limit is as described in
the first example above. If, on the other hand, it approaches the
component $b=0$, whose general member corresponds to a curve $C
\subset \P^1 \times \P^1$ containing the fiber $x=0$, we get a limit
as depicted in the third example. And finally, if $\delta$ approaches
(generically) a point in the intersection of these two components, we
get an example of the third type.

\end{document}